\documentclass[12pt,a4paper]{article}
\usepackage{amssymb,amsmath}
\newtheorem{lemma}{\sl LEMMA}[section]
\newtheorem{corollary}{\sl COROLLARY}[section]
\newtheorem{remark}{\sl REMARK}[section]
\newtheorem{theorem}{\sl THEOREM}[section]
\newtheorem{defini}{\sl DEFINITION}
\newcommand{\No}{\#\ }
\renewcommand{\S}{\mathhexbox{278}\ }
\newcommand{\epr}{\hfill $\square$}
\newcommand{\atp}[2]{\genfrac{}{}{0pt}{}{#1}{#2}}
\newcommand{\R}{\mathbb R}
\newcommand{\N}{\mathbb N}
\newcommand{\Lc}{\mathcal L}
\newcommand{\Kc}{\mathcal K}
\newcommand{\Ic}{\mathcal I}
\newcommand{\Gc}{\mathcal G}
\newcommand{\Oc}{\mathcal O}
\newcommand{\Ec}{\mathcal E}
\newcommand{\Xc}{\mathcal X}
\newcommand{\Sc}{\mathcal S}
\newcommand{\Tc}{\mathcal T}
\newcommand{\lm}{\mathrm {lim}}
\newcommand{\lc}{\mathrm {loc}}
\newcommand{\ms}{\mathrm {mes}}
\newcommand{\spp}{\mathrm {supp}}
\begin{document}
\allowdisplaybreaks
{\small\it
\hfill
TO THE MEMORY OF A.S.KALASHNIKOV }
\bigskip
\begin{center}
{\bf ON THE THEORY OF DISCONTINUOUS SOLUTIONS TO SOME
STRONGLY DEGENERATE PARABOLIC EQUATIONS}
\end{center}
\bigskip
\centerline{\it Yu.G.Rykov}
\smallskip
{\small\it
\begin{center}
Keldysh Institute of Applied Mathematics\\
Miusskaya sq.4, 125047 Moscow Russia\\
E-mail: rykov@spp.keldysh.ru
\end{center}
}
%
%%  ABSTRACT
\bigskip
{\small
\centerline{ABSTRACT}
\smallskip
It is studied the Cauchy problem for the equations of Burgers' type but
with bounded dissipation flux $$u_t+f(u)_x=Q\left(u_x\right)_x,\quad
(t,x)\in \R_+\times \R\ ,$$ where $Q^\prime >0$, $\max |Q(s)|<+\infty$.
Such equation degenerates to hyperbolic one as the velocity gradient tends
to infinity. Thus the discontinuous solutions are permitted. In the paper
the definition of the generalized solution is given and the existence
theorem is established in the classes of functions close to ones of
bounded variation. The main feature of used a priori estimates is the fact
that one needs to estimate only $Q\left(u_x\right)$ which allows to have
in fact {\it arbitrary local growth} of the velocity gradient. The
uniqueness theorem is proven for essentially narrower class of piecewise
smooth functions with regular behavior of discontinuity lines.
}

%
%  SECTION 1
%
\section{Introduction}

It is studied the Cauchy problem to the equations of the following type
\begin{equation}
\Lc u\equiv u_t+f(u)_x-Q\left(u_x\right)_x=0 \label{1}
\end{equation}
in the strip $(t,x)\in\Pi_T\equiv [0,T]\times \R$ with initial
conditions
\begin{equation}
u(0,x)=u_0(x)\ ,\quad x\in \R\ , \label{2}
\end{equation}
where $max_{x\in \R} |u_0(x)|\equiv M<+\infty$.

We will consider the problem of existence of generalized solutions
$u(t,x)\in L_1(\Pi_T)$ to the equation (\ref{1}) from the class
$BV_{C^1}(\R)$ for almost all $t\in [0,T]$; through $BV_{C^1}(\R)$ one
denotes the set of such functions $g(x)\in BV_\lc (\R)\bigcap$ $L_1(\R)$
that $g\in C^1$ on the set of full measure in $\R$. The initial function
$u_0(x)$ satisfies the following conditions: 1) there exists such compact
set $K$ that $u_0(x)\in C^2(K\setminus \{x_i\})$, where $\{x_i\}$ --- the
finite number of points; 2) $u_0^\prime(x_i\pm 0)=0$; 3) $u_0(x)\in W^2_1(
\R\setminus K)$, $Q(u_0^\prime(x))\in W_1^1(\R\setminus K)$. (Thus the
initial function can have discontinuities in finite number of points.) One
will say that $u_0(x)$ belongs to the class $BV_{C^1}^+(\R)$.

We will also consider the uniqueness problem but with respect to
essentially narrower class of functions $\Kc$. Roughly speaking
this is the class of piecewise smooth functions with sufficiently regular
behavior of discontinuity lines (the exact definition will be given in \S
2). Such restriction to the class of functions has most likely the
technical character because of arising difficulties in the proof of
general uniqueness theorem. Our proof is based on the essential fact that
for (\ref{1}) one can define the analog of the concept of characteristic
lines (in case of hyperbolic equations) --- level lines of the function
$u(t,x)$.

Assume that functions $f,Q$ satisfy the following conditions: $f\in
C^1$, $Q\in C^2$, $f(0)=Q(0)=0$, $Q^\prime >0$, $Q(-\infty)\equiv
Q_{-\infty}=const<0$, $Q(+\infty)\equiv Q_{+\infty}=const>0$. Let us also
introduce the notations
$$
\bar Q\equiv \max\left(|Q_{-\infty}|, |Q_{+\infty}|\right),
F(M)\equiv \max_{|s|\le M}|f^\prime(s)|,
Q_1\equiv \max_{s\in \R}Q^\prime (s)\ .
$$

The equation (\ref{1}) actually is the generalization to the case of
bounded dissipation flux $Q$ of known Burgers equation ($f(u)\equiv
u^2/2$, $Q(s)\equiv\varepsilon s$) which was introduced as simplest
turbulence model (see, for example, \cite{Burg_BNondif74}).  Such
generalization arose relatively recently in the problems of nonlinear
diffusion, phase transitions theory and generalization of Navier-Stokes
equations (see, for example, \cite{Rose_PhRev89}-- \cite{Rose_PhRev90}).
The equation (\ref{1}) is the simplest model which describes the
interaction between the nonlinear convective transport and dissipation
process when the dissipation flux is bounded. As $|u_x|\to\infty$ the
equation (\ref{1}) becomes the first order equation $u_t+f(u)_x=0$. So it
is natural to expect the appearance of discontinuities in the generalized
solutions. The fact of emerging of the hyperbolic properties in solutions
to parabolic problems was studied, for example, in \cite{BeDPa_ARMA92} and
called there strong degeneracy.  The degenerate parabolic equations were
intensively studied, see, for example, the review \cite{Kalash_UMN87} and
references therein.  But as a rule only continuous solutions and their
properties were considered. In the paper \cite{VoHud_MaSb69} the developed
theory included discontinuous solutions, but these solutions appeared only
when the term with higher derivatives turned identically to zero for some
range of values of $t,x,u$. So actually there was no relation between the
propagation of discontinuities and viscosity terms.  The peculiarity of
(\ref{1}) follows from the presence of discontinuities and their effective
interaction with viscous terms. This interaction generates some kind of
'boundary layer' in the vicinity of the discontinuity.

Formally it is possible to use the corresponding theory for hyperbolic
equations \cite{Olei_UMNa57}, \cite{Kruz_MatSb70} to study the equation
(\ref{1}). The advantage to formulate the notion of generalized solution
in terms of integral inequality gives at once the validity of entropy
conditions for discontinuities and as a consequence the uniqueness
theorem. But to prove the existence theorem it is necessary to have the
additional smoothness which in general one does not have
\cite{Rykov_Prep99}. From the other hand the presence of viscosity however
degenerate most likely provides the validity of entropy conditions for
discontinuities. So the integral inequality is not necessary and the main
difficulty will be the proof of general uniqueness theorem.

The equation (\ref{1}) and its more complicated variant when the function
$Q$ is non-monotone was studied in \cite{KurgRo_CPAM97},
\cite{KurgLRo_CPAM98}. There were proven the existence and uniqueness
theorems when one has no discontinuities and were shown a number of
numerical calculations by the small viscosity method to illustrate the
qualitative behavior of the solutions. In the recent work
\cite{GoodKurRo_Nonl99} it was proven that in the case of non-monotone as
well as monotone function $Q$ the discontinuities really emerged. There
were also new series of numerical calculations which were based on the
splitting method.

Let us briefly outline the contents of the paper. In \S 2 one formulates
the definition of the notion of generalized solution and some
preliminaries are shown to justify introduced definition. \S 3 is devoted
to the proof of the existence theorem by small viscosity method. To prove
the theorem one needs only a priori estimates for $Q\left(u_x\right)$ and
not for $u_x$ itself. This is the main feature of used estimates and
allows us to have almost {\it arbitrary local growth} of the velocity
gradient. In \S 4 the uniqueness theorem is proved for the functions from
the class $\Kc$ providing that Oleinik's condition $E$ of
\cite{Olei_UMNa59u} is true.

%
%  SECTION 2
%
\section{The formulation of basic results}

%
%  SUBSECTION 2.1
%
\subsection{The definition of generalized solution and main theorems}

Let us first give the definition of generalized solution to the problem
(\ref{1}), (\ref{2}).

%
%  DEFINITION 1
%
\begin{defini}
Bounded measurable function $u(t,x)$ will be called the generalized
solution to the problem (\ref{1}), (\ref{2}) in $\Pi_T$ iff:

1) there exists such set ${\Ec}\subset [0,T]$, $\ms\
\Ec=0$ that as $t\in [0,T]\setminus \Ec$ the function
$u(t,x)\in BV_{C^1}(\R)$ is defined a.~e. in $\R$ and there
exists
$$
\lim_{\atp{h\rightarrow 0}{h\in [0,1]\setminus \Xc (x)}}
Q\left(\frac{u(t,x+h)-u(t,x-h)}{2h}\right)=Q_\lm (t,x)
$$
for every $x\in \R$, where $\ms\ \Xc (x)=0$,
$Q_\lm (t,x)\in BV_\lc (\R)$ and is
continuous with respect to $x$\ ;

2) for an arbitrary function $\varphi\in C_0^\infty (\Pi_T)$ the following
integral identity holds
$$
\int\!\!\int\limits_{\Pi_T}\left\{u(t,x)\varphi_t+f(u(t,x))\varphi_x-
Q_\lm (t,x)\varphi_x\right\}dxdt = 0\ ;
$$

3) for every segment $[a,b]\subset \R$
$$
\lim_{\atp{t\rightarrow 0}{t\in [0,T]\setminus \Ec}}
\int\limits_a^b|u(t,x)-u_0(x)|dx=0\ .
$$
\end{defini}

Define now the functional classes $\Kc_0$ and $\Kc$ which we
will use for the proof of uniqueness theorem.

%
%  DEFINITION 2
%
\begin{defini}
One will say that piecewise smooth function $u(t,x)$ belongs to the class
$\Kc_0$ iff for every $T>0$ as $0<t<T$ the following conditions
hold:
\begin{itemize}
\item[\it i).] At any point of discontinuity except finite number of
points there exist one-sided limits $u(t,x_i(t)\pm
0)\equiv u^\pm(t,x_i(t))$, $(i=1,\dots,N)$, $u^-\ne u^+$. (Thus the
discontinuity lines can intersect only at finite number of points.)
\item[\it ii).] There exists such $\delta>0$ that for every line of
discontinuity or non-smoothness $x_i(t)$, $(i=1,\dots,N_1)$ the equation
$u^\pm(t,x_i(t))=c=const$ has only finite number of solutions for almost
every $|c|\in [0,\delta]$.
\end{itemize}
\end{defini}

%
%  DEFINITION 3
%
\begin{defini}
One will say that the function $u(t,x)$ belongs to the class $\Kc$
iff for every $T>0$ as $0<t<T$ the following conditions hold:
\begin{itemize}
\item[\it i).] The function $u(t,x)\in C^2(\R^2)$ everywhere except
finite number of lines $x_i(t)$, $i=1,\dots ,N$ which themselves belong
to the class $C^2$. Moreover $\sup\limits_{[0,T]}|u(t,R)|\rightarrow
0$ as $|R|\rightarrow 0$.
\item[\it ii).] At every discontinuity point except finite number of them
there exist one-sided limits $u(t,x_i(t)\pm 0)\equiv u^\pm$, $u^-\ne
u^+$.
\item[\it iii).] The function $\widehat Q(t,x)\equiv Q(u_x(t,x))$ as
\hspace*{\fill}\linebreak
$(t,x)\in \R^2\setminus \bigcup\limits_{1\le i\le N,t\in
\R^+}(t,x_i(t))$, otherwise $\widehat Q(t,x)\equiv Q_{-\infty}$ as
$u^->u^+$ and $\widehat Q(t,x)\equiv Q_{+\infty}$ as $u^-<u^+$, is
continuous. Moreover $\sup\limits_{[0,T]}|\widehat Q(t,x)|\rightarrow 0$
as $|R|\rightarrow 0$.
\item[\it iv).] The difference of any two functions from the class
$\Kc$ belongs to the class $\Kc_0$.
\end{itemize}
\end{defini}

%
%  REMARK 2.1
%
\begin{remark}
Let us note that the integral identity 2) and properties of the function
$Q_\lm (t,x)$ from the Definition 1 imply the Hugoniot
condition at the discontinuity lines $y(t)$ for the functions belonging to
the class $\Kc$
$$
\dot y(t)=\frac{f(u^+(t,y(t)))-f(u^-(t,y(t)))}{u^+(t,y(t))-u^-(t,y(t))}\ .
$$
\end{remark}

Remind also necessary for us condition $E$ from the paper
\cite{Olei_UMNa59u}. Let us introduce the notation
$$
l(u)\equiv f(u^-)+(u-u^-)\frac{f(u^+)-f(u^-)}{u^+-u^-}\ .
$$

%
%  DEFINITION 4
%
\begin{defini}
One will say that the generalized solution $u(t,x)$ from the class
$\Kc$ to the problem (\ref{1}), (\ref{2}) satisfies the condition
$E$ iff at every point of discontinuity of the function $u(t,x)$ except
finite number of points the following condition holds: if $u^->u^+$ then
$l(u)\ge f(u)$ for $u\in [u^+,u^-]$; if $u^-<u^+$ then $l(u)\le f(u)$ for
$u\in [u^-,u^+]$.
\end{defini}

The main result of the present paper is the proof of the following
theorems.

%
%  THEOREM 2.1
%
\begin{theorem}
For the Cauchy problem (\ref{1}), (\ref{2}) there exists the generalized
solution $u(t,x)$ in the sense of Definition 1.
\end{theorem}

%
%  THEOREM 2.2
%
\begin{theorem}
If the generalized solution $u(t,x)$ to the Cauchy problem (\ref{1}),
(\ref{2}) belongs to the class $\Kc$ and satisfies the condition
$E$ then it is unique.
\end{theorem}

Let us now study some particular solutions to the equation (\ref{1}) to
demonstrate the validity of Definition 1.

%
%  SUBSECTION 2.2
%
\subsection{The traveling wave solution.}

Suppose that $f^{\prime\prime}>0$ and let seek the solution of (\ref{1})
in the form
\begin{equation}
u(t,x)=b(x-st)\equiv b(\xi)\ ,\label{3}
\end{equation}
where $b(-\infty)=b_-=const$, $b(+\infty)=b_+=const$, $b^\prime\to 0\
\mbox{as}\ |\xi|\to\infty$. After the substitution to (\ref{1}) one has
$$
-sb^\prime(\xi)+f(b(\xi))^\prime-Q(b^\prime (\xi))^\prime =0\ ,
$$
i.e.
$$
Q(b^\prime )=f(b)-f(b_-)-s(b-b_-)\equiv\widehat f(b)\ .
$$
Hence
\begin{equation}
s=\frac{f(b_+)-f(b_-)}{b_+-b_-}\ ;\
\int\limits_{b_0}^b\frac{dB}{Q^{-1}
\left(f(B)-f(b_-)-s(B-b_-)\right)}=\xi\ ,\label{4}
\end{equation}
where $b_0$ is some arbitrary constant.

Suppose $b_->b_+$ and $m\equiv\min_{b\in [b_+,b_-]}\widehat f(b)$.
Then one encounters two cases:

1) If $|m|\le |Q_{-\infty}|$ then formula (\ref{4}) gives the continuous
solution to the equation (\ref{1}) in the form (\ref{3}), for
definiteness it is set that $b_0=(b_++b_-)/2$.

2) If $|m|>|Q_{-\infty}|$ then, denoting through $b_2<b_1$ such values of
$b$ that $\widehat f(b_i)=m$, $i=1,2$, one obtains the discontinuous
solution to the equation (\ref{1}) in the form (\ref{3}) by formula
(\ref{4}):  $b_0=b_1$ as $\xi<0$, $b_0=b_2$ as $\xi>0$. Thus there is the
discontinuity of strength $b_1-b_2$ at point $\xi=0$.

So we discover that in general as $b_--b_+$ is sufficiently large the
discontinuous traveling wave can exist. But formula (\ref{4}) shows that
the function $\widehat Q(\xi)$ : $\widehat Q(\xi)=Q(b^\prime )$ as
$\xi\not =0$, $\widehat Q(\xi)=Q_{-\infty}$ as $\xi =0$, is continuous.
Thus the hyperbolic properties (emerging of discontinuities) and parabolic
properties (smoothness of the graph of the function on the plane
$(\xi,b)$) are combined.

Further, the function $\widehat Q(\xi)$ has the oscillation of the
strength $|Q_{-\infty}|$ in the traveling wave because $Q(b^\prime
(-\infty))=0$, $Q(b^\prime (0))= Q_{-\infty}$, $Q(b^\prime (+\infty))=0$.
Suppose that the typical singularities of (\ref{1}) are similar to the
traveling wave. Than the function $Q(u_x)$ is continuous at the
discontinuity but oscillates. As one has the set of discontinuities one has
the set of oscillations as well.  Suppose at some moment of time $t_0$ we
have the countable number of discontinuities within the finite segment.
Then we have the oscillations of $sin(1/x)$ type and $Q(u_x)$ is now
discontinuous. But at any moment of time $t_0+\Delta\tau$ after
arbitrarily small time interval $\Delta\tau$ these oscillations disappear.
As two discontinuities merge the continuity of $Q(u_x)$ with respect to
$t$ fails because instantaneously single oscillation arises instead of two
ones.

Finally in case $b_-<b_+$ formula (\ref{4}) shows that the solution of
traveling wave type does not exist.

%
%  SUBSECTION 2.3
%
\subsection{The example of the solution with discontinuous initial data}

Consider the equation (\ref{1}) and let $f\equiv 0$. Suppose the function
$Q(s)$ is such that $Q(s)=-Q(-s)$, $Q^\prime (s)\sim K|s|^{-\beta}$ as
$|s|\to\infty$; $\beta>0$, $K>0$. Suppose also that we are given with the
initial step-like function
$$
u(0,x)=\left\{
\begin{aligned}
u\equiv const>0 & \quad \mbox{as}\quad x<0 \\
0 & \quad \mbox{as}\quad x>0
\end{aligned}\right.
$$
and let seek the solution of (\ref{1}) with $f\equiv 0$ in self-similar
form
\begin{equation}
u(t,x)=\left\{
\begin{aligned}
u-\sqrt{t}h(z) & \quad\mbox{as}\quad x<0 \\
\sqrt{t}h(z) & \quad\mbox{as}\quad x>0
\end{aligned}\right.\ ,\label{5}
\end{equation}
where $z=x/\sqrt{t}$. Then for the function $h(z)$ one obtains the
following ordinary differential equation
$$
h(z)-zh^\prime (z)=2\left[Q\left(h^\prime(z)\right)\right]^\prime
$$
with the condition
$$
|h(z)|\to 0\ \mbox{при}\ |z|\to\infty\ .
$$
We are interested in the behavior of the solution in the neighborhood of
the point $z=0$. Representing there the unknown function $h(z)$ in the
form
$$
h(z)=h(0)+Az^\alpha+\ldots\ ,
$$
we obtain that $\alpha=(\beta-2)/(\beta-1)$ as $\beta>2$. It is easy to see
that $0<\alpha<1$ and initial discontinuity is preserved some time
(after the time moment $t^*=\left(u/2h(0)\right)^2$ the
discontinuity disappears and evolution loses self-similar character
(\ref{5})) but its graph has certain smoothness in a sense that
defined above function $\widehat Q$ remains continuous.

From the other hand let us note that generally speaking another
self-similar regime with $\alpha=2$ exists for the same initial data; then
the function $\widehat Q$ is discontinuous at point $x=0$ but this
discontinuity is removable. The analogous situation arises in case
$f\not\equiv 0$ as well. Namely, let us find the solution in Kruzhkov'
sense \cite{Kruz_MatSb70} of the following equation (the function $f$ is
convex)
$$
u_t+f(u)_x=0
$$
with the same as above step-like initial data. We obtain the discontinuous
solution with constant values of $u$ on the left and right of
discontinuity line. It can be easily seen that this solution satisfies
also items 2) and 3) of Definition 1 for the generalized solutions to
equation (\ref{1}). But the solution does not satisfy item 1) because the
corresponding function $\widehat Q$ is discontinuous (although the
discontinuity is yet removable). Thus the restriction of continuity to the
function $Q_\lm$ from item 1) of Definition 1 is necessary to
prove the uniqueness theorem.

%
%  SECTION 3
%
\section{The existence of generalized solutions.}

Let us approximate the equation (\ref{1}) by the uniformly parabolic
equation
\begin{equation}
\Lc_\varepsilon u\equiv u_t+f(u)_x-Q(u_x)_x-\varepsilon u_{xx}=0\
, \label{6}
\end{equation}
where $\varepsilon>0$ can be taken arbitrary small.

Multiplying (\ref{6}) by the function $\varphi\in C_0^\infty (\Pi_T)$
and integrating with respect to $\Pi_T$ one immediately gets as
$\varepsilon\rightarrow 0$ the integral identity {\it 2)} from the
Definition 1 providing $|u|$ is bounded and appropriate convergence as
$\varepsilon\rightarrow 0$ are valid.

Let us now study the problem (\ref{6}), (\ref{2}). Suppose in addition
$f,Q\in C^3$ and $u_0(x)\in C^4(\Omega)$ for every domain $\Omega\Subset
\R$. Then it is known \cite{LadSolUr_BQuas67} that under supposed
requirements there exists unique classical solution $u^\varepsilon(t,x)$
to the problem (\ref{6}), (\ref{2}) and $u^\varepsilon(t,x)\in
C^{4,2}(\bar Q_T)$ for every bounded cylinder $Q_T\subset\Pi_T$. Further
we will get the uniform with respect to $\varepsilon$ estimates for the
solutions to the problem (\ref{6}), (\ref{2}) which do not depend on the
additional smoothness of $f,Q$ and $u_0(x)$. Thus our results will be
valid for the initial data pointed out in (\ref{2}).

Applying maximum principle \cite{LadSolUr_BQuas67} one immediately
obtains
\begin{equation}
|u^\varepsilon(t,x)|\le M\equiv\max_{\R}|u_0|\ .\label{7}
\end{equation}
Consider the linear equation
\begin{equation}
Lz\equiv z_t+\left(a(t,x,\varepsilon)z\right)_x-
\left(b(t,x,\varepsilon)z_x\right)_x=F(t,x)\label{8}
\end{equation}
in $\Pi_T$, where $\varepsilon$ is a parameter, $b\ge 0$, $a,b$ are
bounded and continuous with their first derivatives with respect to $x$ in
every cylinder $Q_T\subset\Pi_T$ (not necessary uniformly with respect to
$\varepsilon$), $F(t,x)\in C(\bar \Pi_T)$.

Let us formulate the following theorem which is direct generalization of
corresponding statement from \cite{Volp_MaSb67}.

%
%  THEOREM 3.1
%
\begin{theorem}
1) Suppose $Q_r\equiv \left\{(t,x)\ :\ |x|<r, 0<t<T\right\}$ and $z(t,x)$
is classical solution to the equation (\ref{8}) in $\Pi_T$. Then for any
$\varphi(t,x)\ge 0$, $\varphi\in C^\infty$, $\varphi$ equals zero outside
of the cylinder $Q_r$, the following equality takes place
\begin{equation}\begin{aligned}
\ & \int\!\!\int\limits_{Q_r}L^*\varphi
|z|dxdt+\int\limits_{\R}\varphi(T,x) |z(T,x)|dx\le \\
\ & \int\limits_{\R}\varphi(0,x)|z(0,x)|dx+
\int\!\!\int\limits_{Q_r}\varphi |F|dxdt\ ,
\end{aligned}\label{9}
\end{equation}
where
\begin{equation}
L^*\varphi=-\left(\varphi_t+a(t,x,\varepsilon)\varphi_x+
\left(b(t,x,\varepsilon)\varphi_x\right)_x\right)\ . \label{10}
\end{equation}

2) If additionally one suppose that $z(0,x)\in L_1(\R)$, $F(t,x)\in
L_1(\Pi_T)$ then
\begin{equation}
\int\limits_{\R}|z(t,x)|dx\le\int\limits_{\R}|z(0,x)|dx+
\int\!\!\int\limits_{\Pi_T}|F(t,x)|dxdt \label{11}
\end{equation}
for every $t\in [0,T]$.
\end{theorem}

%
%  REMARK 3.1
%
\begin{remark}
We formulate this generalization only in one-dimensional case which we
need here. It is clearly seen that the similar assertion is true in
multi-dimensional case as well.
\end{remark}

%
%  THEOREM 3.2
%
\begin{theorem}
Suppose $u^\varepsilon(t,x)$ is classical solution to the problem
(\ref{6}), (\ref{2}) in $\Pi_T$. Then the family of functions
$\left\{u^\varepsilon(t,x)\right\}$ is compact in $L_1(\Pi_T)$.
\end{theorem}

{\sl PROOF.} Differentiating (\ref{6}) with respect to $x$ or with respect
to $t$ and applying Theorem 3.1 one obtains the following estimates
\begin{align}
\int\limits_\R |u^\varepsilon|dx &\le\int\limits_\R |u_0|dx\equiv M_0
\notag \\
\int\limits_\R |u^\varepsilon_x|dx &\le
\int\limits_\R |u_0^\prime|dx\equiv M_1
\label{12} \\
\int\limits_\R |u^\varepsilon_t|dx &\le
\int\limits_\R |\varepsilon
u_0^{\prime\prime}+Q(u_0^\prime)^\prime -f(u_0)^\prime|dx\ , \notag
\end{align}
where 'prime' denotes the differentiation with respect to $x$. Taking into
account the conjecture of sufficient smoothness of initial function
$u_0(x)$ one gets the $L_1$-compactness of the family
$\left\{u^\varepsilon(t,x)\right\}$.

\epr

%
%  COROLLARY 3.1
%
\begin{corollary}
\begin{equation}
\int\limits_\R |Q\left(u^\varepsilon_x\right)_x|dx+
\varepsilon
\int\limits_\R |u^\varepsilon_{xx}|dx\le M_2=const\label{13}
\end{equation}
uniformly with respect to $\varepsilon$.
\end{corollary}

{\sl PROOF.} The result follows from the estimates (\ref{12}), equation
(\ref{6}) and inequality $Q^\prime>0$.

\epr

Further for our convenience we omit superscript $\varepsilon$  in case it
does not influence the clearness of presentation.

Differentiating (\ref{6}) with respect to $x$ one obtains the equation for
$v\equiv u_x$
\begin{equation}
v_t+\left(f^\prime(u)v\right)_x=\left(Q(v)_x+\varepsilon v_x\right)_x\ .
\label{14}
\end{equation}
Denote through $v^c$ the cut-off function for $v$
$$
v^c\equiv\left\{\begin{aligned}
c\ &,\quad v\ge c \\
v\ &,\quad -c\le v\le c \\
-c\ &,\quad v\le -c\ ,
\end{aligned}\right.
$$
and through $Q_c(s)$ the following function
$$
Q_c(s)\equiv\left\{\begin{aligned}
Q(c+1)\ &,\quad s\ge c+1 \\
Q^+(s)\ &,\quad c\le s\le c+1 \\
Q(s)\ &,\quad -c\le s\le c \\
Q^-(s)\ &,\quad -c-1\le s\le -c \\
Q(-c-1)\ &,\quad s\le -c-1 \ ,
\end{aligned}\right.
$$
where $Q^+$ and $Q^-$ are chosen such that $Q_c(s)\in C^2(\R)$ and
$Q_c^\prime(s)>0$ as $-c-1<s<c+1$.

Below the sign of $\int$ denotes the integration with respect to $\R$ and
$\int\!\!\int$ denotes the integration with respect to the whole strip
$\Pi_T$.

%
%  THEOREM 3.3
%
\begin{theorem}
Suppose $v(t,x)\equiv u_x(t,x)$ is the classical solution of the equation
(\ref{14}) in $\Pi_T$. Then
$$
\int\limits_0^T\int\limits_{\R}\left(v^c_{\ x}\right)^2 dxdt\le
K(c,T,M,M_1,\bar Q)
$$
uniformly with respect to $\varepsilon >0$.
\end{theorem}

{\sl PROOF.} Multiply (\ref{14}) on $Q_c(v)\eta(x)$, where
$\eta(x)\ge 0$, $\eta\in \Sc$ (the space of rapidly decreasing at
infinity functions), and integrate with respect to $\R$. Then
integrating by parts one obtains
\begin{align*}
\int v_tQ_c(v)\eta dx-\int f^\prime(u)v
\left[Q_c(v)\eta\right]_xdx= \\
-\int\left(Q(v)+\varepsilon v\right)_x\left[Q_c(v)_x\eta+
Q_c(v)\eta^\prime\right]dx\ .
\end{align*}
Denote through $\widehat Q_c(s)$ the primitive of the function $Q_c(s)$,
i.e. $\widehat Q_c^\prime(s)=Q_c(s)$. Then
\begin{align*}
\frac{d}{dt}\int \widehat Q_c(v)\eta dx-\int f^\prime(u)vQ_c(v)_x\eta
dx-\int f^\prime(u)vQ_c(v)\eta^\prime dx= \\
-\int\left(Q(v)_x+\varepsilon v_x\right)Q_c(v)_x\eta dx-
\int\left(Q(v)_x+\varepsilon v_x\right)Q_c(v)\eta^\prime dx\ .
\end{align*}
Now let us integrate with respect to the segment $[0,T]$
\begin{align*}
&\left.\int \widehat Q_c(v)\eta(x)dx\right|_0^T-
\int\!\!\int f^\prime(u)vQ_c(v)_x\eta dxdt-
\int\!\!\int f^\prime(u)vQ_c(v)\eta^\prime dxdt=
\\ \\
&-\int\!\!\int Q(v)_xQ_c(v)_x\eta dxdt
-\varepsilon\int\!\!\int v_xQ_c(v)_x\eta dxdt-
\\ \\
&\int\!\!\int\left(Q(v)_x+\varepsilon v_x\right)Q_c(v)\eta^\prime dxdt\ ;
\\ \\
&\int\!\!\int Q^\prime(v)Q_c^\prime v_x^2\eta dxdt+
\varepsilon\int\!\!\int Q_c^\prime(v)v_x^2\eta dxdt=
\\ \\
&-\int\!\!\int\left(Q(v)_x+\varepsilon v_x\right)Q_c(v)\eta^\prime dxdt+
\int\!\!\int f^\prime(u)vQ_c(v)\eta^\prime dxdt+
\\ \\
&\int\!\!\int f^\prime(u)vQ_c(v)_x\eta dxdt-
\left.\int\widehat Q_c(v)\eta dx\right|_0^T\ .
\end{align*}
Further, taking into account that $Q_c^\prime\ge 0$; $Q_c^\prime(v)=0$ as
$|v|\ge c+1$, and applying to the third summand on the right hand side the
inequality $2ab\le a^2+b^2$ one has
\begin{align*}
&\iint\limits_{|v|\le c+1}Q^\prime(v)Q_c^\prime(v)v_x^2\eta dxdt
\le\bar Q\iint |Q(v)_x+\varepsilon v_x||\eta^\prime| dxdt+
\\ \\
&\bar QF(M)\iint |v||\eta^\prime| dxdt+
\iint\limits_{|v|\le c+1}\left| f^\prime(u)v
\sqrt{\frac{Q_c^\prime(v)}{Q^\prime(v)}}\sqrt{\eta}\cdot
v_x\sqrt{Q_c^\prime(v)Q^\prime(v)}\sqrt{\eta}\right| dxdt+
\\ \\
&\bar Q\int |v|(T)\eta dx+\bar Q\int |v|(0)\eta dx\le
\bar Q\iint |Q(v)_x+\varepsilon v_x||\eta^\prime| dxdt+
\\ \\
&\bar QF(M)\iint |v||\eta^\prime| dxdt+\frac{1}{2}
\iint\limits_{|v|\le c+1}f^\prime(u)^2v^2
\frac{Q_c^\prime(v)}{Q^\prime(v)}\eta dxdt+
\\ \\
&\frac{1}{2}\iint\limits_{|v|\le c+1}Q_c^\prime Q^\prime(v)v_x^2
\eta dxdt+2\bar Q\int |u_0^\prime|\eta dx\ .
\end{align*}
Now take $\eta(x)=exp(-\lambda\sqrt{1+x^2})$ and using the inequality
$|\eta^\prime|\le\lambda\eta$ one gets when $\lambda$ tends to zero
(accounting of (\ref{12}), (\ref{13}))
$$
\frac{1}{2}\iint\limits_{|v|\le c+1} Q^\prime(v)Q_c^\prime(v)v_x^2
dxdt\le\frac{F(M)^2}{2}K_0(c)\iint |v|dxdt+
2\bar Q\int |u_0^\prime|dx\ ,
$$
where
$$
K_0(c)\equiv\max_{|v|\le c+1}\frac{vQ_c^\prime(v)}{Q^\prime(v)}\ .
$$
Further,
\begin{align*}
&\frac{1}{2}\left[\min_{|v|\le c} {Q^\prime}^2(v)\right]
\iint\left(v_{\ x}^c\right)^2 dxdt\le\frac{1}{2}
\iint\limits_{|v|\le c}Q^\prime(v)Q_c^\prime(v)v_x^2 dxdt\le
\\ \\
&\frac{1}{2}
\iint\limits_{|v|\le c+1}Q^\prime(v)Q_c^\prime(v)v_x^2 dxdt\le
\left(\frac{F(M)^2}{2}K_0(c)T+2\bar Q\right)\int |u_0^\prime| dx\ .
\end{align*}

\epr

%
%  THEOREM 3.4
%
\begin{theorem}
For every $\eta\in C_0^\infty(\R)$ the sequence of the functions
\linebreak $\int Q\left(u_x^\varepsilon(t,x)\right) \eta(x)dx$ is compact
in $L_1([0,T])$.
\end{theorem}

{\sl PROOF.} Let us prove uniform boundedness and equicontinuity in
$L_1([0,T])$ of the sequence of functions mentioned above:
\begin{align*}
&i).\quad
\int\limits_0^T\left|\int Q(u_x^\varepsilon)\eta(x)dx\right|dt=
\\ \\
&\int\limits_0^T\left|\int dx\int\limits_0^1
Q^\prime(\theta u_x^\varepsilon)d\theta u_x^\varepsilon\eta(x)\right|dt\le
\int\limits_0^T\!\!\int\left|\int\limits_0^1
Q^\prime(\theta u_x^\varepsilon)d\theta\right| |u_x^\varepsilon||\eta|
dxdt\le \\ \\
&\max_{\R} Q^\prime\cdot T\cdot\max
|\eta|\cdot\int|u_x^\varepsilon|dx\le Q_1\cdot T\cdot\max |\eta|\cdot\int
|u_0^\prime|dx\ ;
\end{align*}

\begin{align*}
&ii).\quad I\equiv\int\limits_0^T\left|\int\left[
Q\left(u_x^\varepsilon(t+\Delta t,x)\right)-
Q\left(u_x^\varepsilon(t,x)\right)\right]\eta(x)dx\right|\le
\\ \\
&\int\limits_0^T\biggl\{\left|\int\left[
Q\left(u_x^\varepsilon(t+\Delta t,x)\right)-
Q_c\left(u_x^\varepsilon(t+\Delta t,x)\right)\right]\eta(x)dx\right|+
\\ \\
&\left|\int\left[
Q_c\left(u_x^\varepsilon(t+\Delta t,x)\right)-
Q_c\left(u_x^\varepsilon(t,x)\right)\right]\eta(x)dx\right|+
\\ \\
&\left|\int\left[
Q_c\left(u_x^\varepsilon(t,x)\right)-
Q\left(u_x^\varepsilon(t,x)\right)\right]\eta(x)dx\right|\biggr\}dt\ .
\end{align*}
For an arbitrary small $\Delta >0$ choose sufficiently large $c>0$ such
that
$$
|Q(v)-Q_c(v)|\le\frac{\Delta}{3T\int |\eta|dx}\ .
$$
This is possible because $Q$ and $Q_c$ tend to the constants at infinity.

Further using the equation (\ref{14}) one has
\begin{align*}
&I\le\frac{2\Delta}{3}+ \int\limits_0^T\left|\int dx\eta(x)
\int\limits_{t}^{t+\Delta t}Q_c\left(u_x^\varepsilon\right)_\tau d\tau
\right|dt=
\\ \\
&\frac{2\Delta}{3}+ \int\limits_0^T\left|
\int\limits_{t}^{t+\Delta t}d\tau\int Q_c^\prime(v)\left[
\left(Q(v)_x+\varepsilon v_x\right)_x-\left(f^\prime (u)v\right)_x\right]
\eta(x)dx\right|dt=
\\ \\
&\frac{2\Delta}{3}+ \int\limits_0^T dt\biggl|
\int\limits_{t}^{t+\Delta t}d\tau\biggl\{\int Q_c^\prime(v)_x\left[
Q(v)_x+\varepsilon v_x-f^\prime (u)v\right]
\eta(x)dx +
\\ \\
&\int\limits_{t}^{t+\Delta t}d\tau\biggl\{ \int Q_c^\prime(v)\left[
Q(v)_x+\varepsilon v_x-f^\prime (u)v\right]
\eta^\prime(x)dx\biggr\}\biggr|\le
\\ \\
&\frac{2\Delta}{3}+ \int\limits_0^T dt
\int\limits_{t}^{t+\Delta t}d\tau\int dx|\eta(x)|\biggl\{
\left|Q_c^{\prime\prime}(v)\left(
Q^\prime(v)+\varepsilon\right)v^2_x\right|+
\\ \\
&\left|Q_c^{\prime\prime}(v)v_xf^\prime(u)v\right|\biggr\}+
\max |\eta^\prime|\cdot\max Q_c^\prime\cdot
\int\left[|Q(v)_x+\varepsilon v_x|+F(M)|v|\right]dx\cdot T\Delta t.
\end{align*}
Hence making the change of variables $\tau^\prime=\tau$,
$t^\prime=-t+\tau$ and using Theorem 3.3 we get
\begin{align*}
&I\le\frac{2\Delta}{3}+ \int\limits_0^{\Delta t}dt^\prime
\int\limits_{t^\prime}^{T+t^\prime}d\tau^\prime\int\limits_{|v|\le c+1}dx
|\eta(x)|\biggl\{\left|
Q_c^{\prime\prime}(v)\left(Q^\prime(v)+\varepsilon\right)v_x^2\right|+
\\ \\
&\left| Q_c^{\prime\prime}(v)v_x f^\prime(u)v\right|\biggr\}+
\Delta t\cdot T\cdot\max |\eta^\prime|\cdot\max Q_c^\prime\cdot
\left(M_2+F(M)M_1\right)\le
\\ \\
&\frac{2\Delta}{3}+ \int\limits_0^{\Delta t}dt^\prime
\int\limits_{0}^{T+1}d\tau^\prime\biggl\{\left[\max_{|v|\le c+1}
Q_c^{\prime\prime}(v)(Q^\prime(v)+1)\right]\int\limits_{|v|\le c+1}dx
|\eta(x)|\left(v_x^{c+1}\right)^2+
\\ \\
&F(M)\left[\max_{|v|\le
c+1}Q_c^{\prime\prime}(v)v\right]\int\limits_{|v|\le c+1}dx
|\eta(x)||v_x^{c+1}|\biggr\}+
\Delta t\cdot\max |\eta^\prime|\cdot const_1\le
\\ \\
&\frac{2\Delta}{3}+ \Delta t\cdot\left[
\max_{|v|\le c+1} Q_c^{\prime\prime}(v)(Q^\prime(v)+1)\right]\cdot
\max |\eta|\cdot K(c+1,T+1,M,\bar Q,M_1)+
\\ \\
&\Delta t\cdot F(M)\cdot\left[
\max_{|v|\le c+1} Q_c^{\prime\prime}(v)v\right]\frac{1}{2}
\left[K(c+1,T+1,M,\bar Q,M_1)+\int\eta^2(x)dx\right]+
\\ \\
&\Delta t\cdot\max |\eta^\prime|\cdot const_1\le
\frac{2\Delta}{3}+\Delta t\cdot const_2\ .
\end{align*}
Now as $c$ is fixed let choose $\Delta t$ in such a way that $\Delta
t\cdot const_2\le\Delta /3$, i.e. $I\le\Delta$.

\epr

%
%  LEMMA 3.1
%
\begin{lemma}
There exists such countable family of functions $\eta_n\in
C_0^\infty(\R),$ that for every $\eta\in C_0^\infty(\R)$ and every
$\delta>0$ there exists such $\eta_n$ that
$\sup_{\R}|\eta_n-\eta|\le\delta$.
\end{lemma}

{\sl PROOF.} Consider the family of all polynomials $P_k(x)$, $x\in\R$
with rational coefficients. Suppose
$$
P_{kl}=\left\{\begin{aligned}
0\quad &,\quad |x|>l,l\in N \\
P_k(x)\quad &,\quad |x|<l,l\in N\ ,
\end{aligned}\right.
$$
and for $m\in \N$
\begin{equation}
P_{klm}=\frac{1}{h}\int\limits_\R \omega\left(\frac{x-y}{h}\right)
P_{kl}(y)dy\ ,\quad h=\frac{1}{m}\ ,\label{15}
\end{equation}
where $\omega\in C_0^\infty$, $\omega\ge 0$, $\int_\R \omega(y)dy=1$. The
family $\left\{P_{klm}(x)\right\}\subset C_0^\infty$ and obviously
countable. Consider some arbitrary function $\eta(x)\in C_0^\infty(\R)$,
suppose $\spp\ \eta(x)\in [-l_1,l_1],l_1\in \N$. Then for every $\delta
>0$ for the segment $[-l_1,l_1]$ by the aid of Weierstrass theorem there
exists the polynomial with rational coefficients $P_{k_1}(x)$ such that
$\sup_{[-l_1,l_1]} |P_{k_1}(x)-\eta(x)|\le\delta /5$. Let us extend
$P_{k_1}(x)$ by zero outside the interval $(-l_1,l_1)$ and denote obtained
function through $P_{k_1l_1}(x)$. Then there exists such $m_1\in \N$ that
$|P_{k_1l_1m_1}(x)- P_{k_1l_1}(x)|\le 4\delta/5$, where $P_{k_1l_1m_1}(x)$
is defined with respect to the formula (\ref{15}). Hence
$|P_{k_1l_1m_1}(x)-\eta(x)|\le\delta$.

\epr

%
%  THEOREM 3.5
%
\begin{theorem}
There exists such subsequence $\{\varepsilon_k\}$ that the sequence
\linebreak
$Q\left(u_x^{\varepsilon_k}(t,x)\right)$ converges in $L^1(\Pi_T)$ to the
function $Q_\lm (t,x)$. Moreover for almost all $t_*\in [0,T]$ the
sequence $Q\left(u_x^{\varepsilon_k}(t_*,x)\right)$ converges in $L^1(\R)$
to the function $Q_\lm (t_*,x)\in L^1(\R)\bigcap BV_\lc (\R)$.
\end{theorem}

{\it PROOF.} Let us take an arbitrary function $\eta_n$ from Lemma 3.1,
then according to the Theorem 3.4 it is possible to choose such
subsequence $\{\varepsilon_k^{(n)}\}$ that there exists the set $\Ec_n$,
$\ms\ \Ec_n=0$ such that $\int Q\left(u_x^{\varepsilon_k^{(n)}}\right)
\eta_n(x)dx$ converges for every $t\in [0,T]\setminus \Ec_n$. With the
help of diagonal process it is possible to choose required subsequence
$\varepsilon_k$ in such a way that $\int Q\left(u_x^{\varepsilon_k}\right)
\eta_n(x)dx$ converges for every $n$ and all $t\in
[0,T]\setminus\bigcup\limits_n \Ec_n$.

Further for an arbitrary $\eta(x)\in C_0^\infty(\R)$ by Lemma 3.1 choose
such $n$ that $\int |\eta_n-\eta|dx\le\Delta/(4\bar Q)$, then
\begin{align*}
&I\equiv\left|\int Q\left(u_x^{\varepsilon_{k_1}}\right)\eta(x)dx-
\int Q\left(u_x^{\varepsilon_{k_2}}\right)\eta(x)dx\right|\le
\\ \\
&\left|\int Q\left(u_x^{\varepsilon_{k_1}}\right)\left(\eta(x)-
\eta_n(x)\right)dx\right|+
\left|\int\left( Q\left(u_x^{\varepsilon_{k_1}}\right)-
Q\left(u_x^{\varepsilon_{k_2}}\right)\right)\eta_n(x)dx\right|+
\\ \\
&\left|\int Q\left(u_x^{\varepsilon_{k_2}}\right)\left(\eta_n(x)-
\eta(x)\right)dx\right|\le 2\bar Q\int |\eta_n-\eta|dx+
\\ \\
&\left|\int\left( Q\left(u_x^{\varepsilon_{k_1}}\right)-
Q\left(u_x^{\varepsilon_{k_2}}\right)\right)\eta_n(x)dx\right|\ .
\end{align*}
Now for fixed $n$ let choose such $\varepsilon >0$ that as
$\varepsilon_{k_1}<\varepsilon$, $\varepsilon_{k_2}<\varepsilon$ the
second integral in the last inequality does not exceed $\Delta /2$. So the
sequence $Q\left(u_x^{\varepsilon_k}\right)$ converges weakly for almost
all $t\in [0,T]$ and therefore it converges in $L_1$ because of estimate
(\ref{13}). Taking into account the boundedness of the function $Q$, by
Lebesgue theorem one obtains the convergence in $L^1(\Pi_T)$. Taking into
account (\ref{13}) it is clear that $\lim_{\varepsilon_k\rightarrow 0}
Q\left(u_x^{\varepsilon_k}\right)\in BV_\lc (\R)$.

\epr

%
%  LEMMA 3.2
%
\begin{lemma}
Suppose one has the sequence of integrable functions $a_n(t)\ge 0$ and
$\int\limits_0^T a_n(t)dt\le C(T)=const$. Then
$A_k(t)\equiv\inf\limits_{n\ge k} a_n(t)<+\infty$ for almost all $t\in
[0,T]$.
\end{lemma}

{\sl PROOF.} One has $A_k(t)\le a_k(t)$, $\int\limits_0^T A_k(t)dt \le
C(T)$, $A_k(t)\le A_{k+1}(t)$. By the B.~Levi theorem $A_k(t)<+\infty$ for
almost all $t\in [0,T]$.

\epr

%
%  LEMMA 3.3
%
\begin{lemma}
Consider an arbitrary sequence of numbers $c_n\rightarrow +\infty$.  Then
for almost all $t\in [0,T]$ there exists such subsequence
$\varepsilon_p(t)\rightarrow 0$ that
\begin{equation}
\int\limits_\R
\left(v_x^{c_n,\varepsilon_p}(t,x)\right)^2dx<+\infty \label{16}
\end{equation}
uniformly with respect to $\varepsilon_p(t)$.
\end{lemma}

{\sl PROOF.} Consider an arbitrary sequence $c_n\rightarrow +\infty$.
According to the theorem 3.3 one has
$$
\int\limits_0^T\!\!\int\limits_\R
\left(v_x^{c_n,\varepsilon}\right)^2dxdt<K(\cdot,c_n)\ ,
$$
where through $(\cdot)$ one denotes the dependence on arguments we do not
worry about at the moment.  By Lemma 3.2 there exist the sets $\Ec_n$,
$\ms\ \Ec_n=0$ such that for every $t\in [0,T]\setminus \Ec_n$
\begin{equation}
V_\alpha(t)\equiv\inf\limits_{0<\varepsilon<\alpha}\int\limits_\R
\left(v_x^{c_n,\varepsilon}\right)^2dx<+\infty\ .\label{17}
\end{equation}
Then for $t\in \Tc\equiv [0,T]\setminus\bigcup\limits_n \Ec_n$ with the
aid of diagonal process one obtains (\ref{16}) from (\ref{17}).

\epr

%
%  THEOREM 3.6
%
\begin{theorem}
The function $Q_\lm (t,x)$ is continuous with respect to $x$ for almost
all $t$.
\end{theorem}

{\sl PROOF.} According to Lemma 3.3 for almost all $t$ there exists such
subsequence $\varepsilon_p(t)$ that the estimate (\ref{16}) is valid.  For
such $t$ let us check the equicontinuity with respect to $x$ (dropping in
the notations index $p$ and argument $t$)
\begin{align*}
&I\equiv\sup_{x\in K\Subset \R}\left|
Q\left(u_x^{\varepsilon}(t,x+\Delta x)\right)-
Q\left(u_x^{\varepsilon}(t,x)\right)\right|\le
\\ \\
&\sup_{x\in K\Subset \R}\left|
Q\left(u_x^{\varepsilon}(t,x+\Delta x)\right)-
Q_c\left(u_x^{\varepsilon}(t,x+\Delta x)\right)\right|+
\\ \\
&\sup_{x\in K\Subset \R}\left|
Q_c\left(u_x^{\varepsilon}(t,x+\Delta x)\right)-
Q_c\left(u_x^{\varepsilon}(t,x)\right)\right|+
\\ \\
&\sup_{x\in K\Subset \R}\left|
Q_c\left(u_x^{\varepsilon}(t,x)\right)-
Q\left(u_x^{\varepsilon}(t,x)\right)\right|\equiv I_1+I_2+I_3\ .
\end{align*}
Choose $c=c_n$ in such a way that $I_1+I_3\le 2\Delta /3$; and by H\"older
inequality
\begin{align*}
&I_2=\sup_{x\in K\Subset \R}\left|
\int\limits_x^{x+\Delta x}\frac{\partial}{\partial\xi}
Q_c\left(u_\xi^{\varepsilon}(t,\xi)\right)d\xi\right|=
\sup_{x\in K\Subset \R}\left|
\int\limits_x^{x+\Delta
x}Q_c^\prime\left(u_\xi^{\varepsilon}(t,\xi)\right) v_\xi^\varepsilon
d\xi\right|\le
\\ \\
&Q_1\sup_{x\in K\Subset \R}
\int\limits_x^{x+\Delta x}\left| v_\xi^{c,\varepsilon}\right| d\xi\le
Q_1\Delta x^{1/2}\sup_{x\in K\Subset \R}\left(
\int\limits_{\R}\left| v_\xi^{c,\varepsilon}\right|^2
d\xi\right)^{1/2}=
\\ \\
&Q_1\Delta x^{1/2}\left(
\int\limits_{\R}\left| v_x^{c,\varepsilon}\right|^2
dx\right)^{1/2}\le const(c_n)\cdot\Delta x^{1/2}\ .
\end{align*}
Now taking $\Delta x$ small enough to provide that right hand side is less
than $\Delta /3$ we get necessary equicontinuity and assertion of Theorem
3.6.

\epr

%
%  THEOREM 3.7
%
\begin{theorem}
For almost all $t\in [0,T]$
\begin{equation}
Q_\lm (t,x)=\lim_{h\rightarrow 0}
Q\left(\frac{u(t,x+h)-u(t,x-h)}{2h}\right)\ .\label{18}
\end{equation}
\end{theorem}

{\sl PROOF.} It follows from the proof of the Theorem 3.6 that for almost
all $t\in [0,T]$ there exists uniformly converging (for each $t$ its own)
subsequence. Let us fix the value of $t\in [0,T]$ and denote such
subsequence again through $\{u^\varepsilon(t,x)\}$.

Suppose $S_t\equiv\left\{x\in \R\ :\ |u_x^\varepsilon(t,x)|\rightarrow
+\infty\right\}$. Since $Q\left(u_x^\varepsilon\right)\rightarrow Q_\lm
(t,x)$ and the function $Q$ is monotone one encounters two cases:

$i).\quad x_0\in \R\setminus S_t$

Since the function $Q_\lm (t,x)$ is continuous there exists such interval
$(a,b)\ni x_0$ that $Q\left(u_x^\varepsilon\right)\rightarrow Q_\lm (t,x)$
and $Q_-<Q_\lm (t,x)<Q_+$. Then for $(a,b)$ the derivative $Q^\prime\ge
const>0$ and $\{u_x^\varepsilon\}$ converges uniformly in $K\Subset
(a,b)$, the function $u=\lim u^\varepsilon$ has continuous derivative and
$Q\left(u_x^\varepsilon\right)\rightarrow Q(u_x)$.

$ii).\quad x_0\in S_t\ ,\ |u_x^\varepsilon(t,x_0)|\rightarrow +\infty$

As it can be seen from item $i)$ the set $S_t$ is closed. Suppose there
exists the interval $(a,b)\ni x_0$ such that $|u_x^\varepsilon(t,x)|\to
+\infty$ for some set of positive measure in $(a,b)$. As far as $Q_\lm
(t,x)$ is continuous one can think for definiteness that
$u_x^\varepsilon(t,x)\rightarrow +\infty$ uniformly with respect to
$\varepsilon$. Further
$$
u^\varepsilon(t,x)=u^\varepsilon(t,x_0)+
\int\limits_{x_0}^xu_x^\varepsilon(t,x)dx\ ,
$$
i.e. $|u^\varepsilon|\rightarrow +\infty$ for the set of positive measure
in some segment $[a_1,b_1]\subset (a,b)$, $x_0\not\in [a_1,b_1]$. But this
contradicts the convergence of $\{u^\varepsilon\}$ in the space $L^1$.
Thus we have proved that $\ms\ S_t=0$.

$iii).$ Suppose for definiteness $u_x^\varepsilon(t,x_0)\rightarrow
+\infty$.

Because of continuity of the function $Q_\lm (t,x)$ for all $x$ from some
neighborhood $U_N(x_0)$ the estimate
$u_x^\varepsilon(t,x)>N$ is true uniformly with respect to
$\varepsilon>0$, where $N>0$ is an arbitrary large number. Then
$$
u^\varepsilon(t,x_0+h)-u^\varepsilon(t,x_0-h)=
\int\limits_{x_0-h}^{x_0+h}u_\xi^\varepsilon(t,\xi)d\xi\ge N\cdot 2h
$$
uniformly with respect to $\varepsilon>0$. Tending $\varepsilon$ to $0$
one has
$$
\frac{u(t,x_0+h)-u(t,x_0-h)}{2h}\ge N
$$
for almost all sufficiently small $h>0$.

\epr

To extend our results up to the initial data (\ref{2}) it is enough to
approximate $u_0\in BV_{C^1}^+$ by smooth initial functions and obtain
uniform estimates for the integrals from (\ref{12}). So consider some
function $\omega(z)\ge 0$, $\omega\in %C_0^\infty$, $supp\ \omega\subset
[-1,1]$, $\int\omega(z)dz=1$ and averaged functions
$$
u_0^h(x)=\frac{1}{h}\int\omega\left(\frac{x-y}{h}\right)u_0(y)dy\ .
$$

%
% LEMMA 3.4
%
\begin{lemma}
Uniformly with respect to $h$ one has the following estimates
\begin{align}
&i).\quad \int\left|\left(u_0^h\right)_x\right|dx\le const_1
\notag \\ \notag \\
&ii).\quad \int\left|Q\left(\left(u_0^h\right)_x\right)_x\right|dx\le
const_2 \label{19} \\
\notag \\
&iii).\quad \int\left|\left(u_0^h\right)_{xx}\right|dx\le const_3/h
\notag
\end{align}
\end{lemma}

{\sl PROOF.} In view of properties of the function $u_0(x)\in BV_{C^1}^+$
it is enough to check inequalities (\ref{19}) operating within small
vicinity of the set of discontinuity points $\{x_i\}$. In addition without
loss of generality it can be considered that $\{x_i\}$ consists of single
point $x_0$.

Let introduce the notation $[u]\equiv |u_0(x_0+0)-u_0(x_0-0)|$. Further
\begin{align*}
&i)\quad \int\limits_{x_0-\delta}^{x_0+\delta}
\left|\left(u_0^h\right)_x\right|dx=
\int\limits_{x_0-\delta}^{x_0+\delta}\left|\frac{1}{h}
\int\limits_{x-h}^{x+h}\frac{\partial}{\partial y}
\omega\left(\frac{x-y}{h}\right)u_0(y)dy\right|dx=
\\ \\
&\int\limits_{x_0-\delta}^{x_0+\delta}\left|\frac{1}{h}
\int\limits_{x-h}^{x_0}\frac{\partial}{\partial y}
\omega\left(\frac{x-y}{h}\right)u_0(y)dy+\frac{1}{h}
\int\limits_{x_0}^{x+h}\frac{\partial}{\partial y}
\omega\left(\frac{x-y}{h}\right)u_0(y)dy\right|dx=
\\ \\
&\int\limits_{x_0-\delta}^{x_0+\delta}\Biggl|\frac{1}{h}
\omega\left(\frac{x-x_0}{h}\right)\left(u_0(x_0-0)-u_0(x_0+0)\right)-
\\ \\
&\frac{1}{h}\int\limits_{x-h}^{x+h}
\omega\left(\frac{x-y}{h}\right)u_0^\prime(y)dy\Biggr|dx\le
[u]+2\delta\max\limits_{
\atp{x_0-\delta\le x<x_0}{x_0< x\le x_0+\delta}}
\left|u_0^\prime(x)\right|\ .
\end{align*}
\begin{align*}
&ii)\quad \Ic\equiv\int\limits_{x_0-\delta}^{x_0+\delta}
\left|Q^\prime\left(\left(u_0^h\right)_x\right)\left(u_0^h\right)_{xx}
\right|dx= \\ \\
&\int\limits_{x_0-\delta}^{x_0+\delta}Q^\prime\left(\frac{1}{h}
\int\limits_{x-h}^{x+h}\frac{\partial}{\partial x}
\omega\left(\frac{x-y}{h}\right)u_0(y)dy\right)\cdot
\left|\frac{1}{h}
\int\limits_{x-h}^{x+h}\frac{\partial^2}{\partial y^2}
\omega\left(\frac{x-y}{h}\right)u_0(y)dy\right|dx\ .
\end{align*}
Estimate the second derivative of averaged function
\begin{align*}
&\left|\left(u_0^h\right)_{xx}\right|=\left|\frac{1}{h}
\int\limits_{x-h}^{x+h}\frac{\partial^2}{\partial y^2}
\omega\left(\frac{x-y}{h}\right)u_0(y)dy\right|= \\ \\
&\Biggl|\frac{1}{h^2}\omega^\prime
\left(\frac{x-x_0}{h}\right)\left(u_0(x_0-0)-u_0(x_0+0)\right)+ \\ \\
&\frac{1}{h}\int\limits_{x-h}^{x_0}
\omega\left(\frac{x-y}{h}\right)u_0^{\prime\prime}(y)dy+
\frac{1}{h}\int\limits_{x_0}^{x+h}
\omega\left(\frac{x-y}{h}\right)u_0^{\prime\prime}(y)dy\Biggr|\le \\ \\
&\frac{1}{h^2}|\omega^\prime|
\left(\frac{x-x_0}{h}\right)[u]+
\max\limits_{
\atp{x_0-\delta\le x<x_0}{x_0< x\le x_0+\delta}}
\left|u_0^{\prime\prime}(x)\right|\ ,
\end{align*}
taking into account the vanishing of one-sided derivatives of the function
$u_0$ at discontinuity points. From estimates of item $i)$ one also infers
\begin{align*}
&\left|\left(u_0^h\right)_{x}-\frac{1}{h}\omega
\left(\frac{x-x_0}{h}\right)[u]\right|\le
\max\limits_{
\atp{x_0-\delta\le x<x_0}{x_0< x\le x_0+\delta}}
\left|u_0^{\prime}(x)\right|\ .
\end{align*}
Therefore
\begin{align*}
&\left|Q^\prime\left(\left(u_0^h\right)_{x}\right)-
Q^\prime\left(\frac{1}{h}\omega\left(\frac{x-x_0}{h}\right)[u]
\right)\right|\le
const_1\cdot\max\limits_{
\atp{x_0-\delta\le x<x_0}{x_0< x\le x_0+\delta}}
\left|u_0^{\prime}(x)\right|\le const_2\cdot h\ .
\end{align*}
Hence
\begin{align*}
&\Ic\le const\left([u]+h+\delta\right)+
\\ \\
&\left|\int\limits_{x_0-\delta}^{x_0+\delta}
Q^\prime\left(\frac{1}{h}\omega\left(\frac{x-x_0}{h}\right)[u]\right)
\left(\frac{1}{h^2}\omega^\prime\left(\frac{x-x_0}{h}\right)[u]
\right)dx\right|\le \\ \\
&const\left([u]+h+\delta\right)+2\bar Q\ .
\end{align*}

$iii).\quad$ Taking into account the estimates from item $ii)$ one has
\begin{align*}
&\int\limits_{x_0-\delta}^{x_0+\delta}
\left|\left(u_0^h\right)_{xx}\right|dx\le
\frac{const_1}{h}[u]+const_2\cdot\delta\ .
\end{align*}

Now if $\varepsilon$ and $\delta$ are taken of order $h$ then one obtains
the estimates (\ref{19}) and the integrals on the right hand side of
(\ref{12}) can be estimated independently on $h$.

\epr

%
%  REMARK 3.2
%
\begin{remark}
We can see from the proved Lemma that the number of discontinuities in
general should be finite because each discontinuity however its strength
put the contribution $2\bar Q$ to the variation of the function
$Q(u_0^\prime)$. If nevertheless one allows to exist infinite number of
discontinuity points than probably it is reasonable to require that the
corresponding function $\widehat Q(u_0^\prime)$ should be continuous (see
\S 2.2), i.e. it should have infinite one-sided derivatives.
\end{remark}

Thus the Theorem 2.1 --- the existence theorem --- has been proved.

%
%  SECTION 4
%
\section{On the uniqueness of generalized solutions.}

Now let us prove the uniqueness Theorem 2.2.

Suppose there exist two solutions $u(t,x)$, $v(t,x)$ from the class
$\Kc$ to the problem (\ref{1}), (\ref{2}) in the sense of
Definition 1. Consider the difference $u-v\equiv w$. Then $w\in
\Kc_0$ and
\begin{equation}
\iint\limits_{\Pi_T}\left\{
(u-v)\varphi_t+\left[f(u)-f(v)\right]\varphi_x+
\left[Q_\lm^u-Q_\lm^v\right]\varphi_x
\right\}dxdt=0\ ,\label{20}
\end{equation}
here $\varphi\in C_0^\infty(\Pi_T)$ and $Q^u_\lm (t,x)$
means the function $Q_\lm$ from Definition 1 which corresponds
to the function $u(t,x)$.

Because of the properties $i)$ in Definition 2 of the functions from the
class $\Kc_0$ the discontinuity lines of the function $w$ can
intersect only at finite number of points. Hence the strip $\Pi_T$ can be
decomposed in finite number of nonoverlapping domains $\Oc_i$,
$i=1,\dots,m$ where the function $w\in C^2$. Consider the level lines
$w=\alpha=const$, $\alpha\in [0,\delta]$, $\delta$ is sufficiently small.
With the help of Sard theorem (see, for example, \cite{PostBLec87}) one
infers that for almost every $\alpha\in [0,\delta]$ in each domain
$\Oc_i$ the level lines of the function $w$ consist of finite
number of regular curves. In consequence of the property $ii)$ of
Definition 2 for almost every $\alpha$ the level lines and discontinuity
lines intersect only at finite number of points.

It follows from (\ref{20}) and item $iii)$ of Definition 3 that for the
discontinuities Hugoniot conditions hold and for every piecewise
$C^1$-contour $\Gamma$ in the strip $\Pi_T$ the following integral
equality holds (the orientation of the plane has been chosen as $(x,t)$)
\begin{equation}
\oint\limits_\Gamma\left(\left[f(u)-f(v)\right]-\left[
Q_\lm^u-Q_\lm^v\right]
\right)dt-(u-v)dx=0\ .\label{21}
\end{equation}

Consider the connected components $\Gc_i$, $i=1,\dots,m_1$ of the domain
$w>\alpha$ (for the connected components $\Gc_i$, $i=m_1+1,\dots,m_1+m_2$
of the domain $w<\alpha$ all considerations are similar). In view of the
decrease at infinity of the functions $u$, $v$, $Q_\lm^u$, $Q_\lm^v$ (see
properties i), iii) of the Definition 3) one can assume that the domain
$w>\alpha$ is bounded.  Indeed the integral of type (\ref{21}) with
respect to the straight lines $|x|=R$ tends to zero as $R\rightarrow 0$.
Consider any connected component $\Gc_i$, without loss of generality one
can take $\Gc_1$, such that $\partial\bar \Gc_1\bigcap
\{t=T\}\not=\varnothing$. So the boundary of the domain $\Gc_1$ will
consist of four type of curves:
\begin{itemize}
\item[\it a).]
The segments $I_k$, $k=1,\dots,l$ of the straight line $t=T$.
\item[\it b).]
The segments of the straight line $t=0$.
\item[\it c).]
Level lines $w=\alpha$, $w$ is continuous.
\item[\it d).]
The lines of discontinuities of the function $w$, where
$w$ crosses the value $\alpha$ step-wise. (Here without loss of generality
it can be reckoned that only function $u$ has the discontinuities.)
\end{itemize}
Apply the formula (\ref{21}) to the contour $\Gamma=\partial\bar
\Gc_1$. Then observe the following.

{\it i).} The integral with respect to lines of type a). gives
$\sum\limits_k\int\limits_{I_k}(u-v)dx$.

{\it ii).} The integral with respect to lines of type b). gives 0 because
of the coincidence of initial values for the functions $u$ and $v$.

{\it iii).} Consider the integral with respect to the line $x_{iii}$ of
type c). Without loss of generality it can be assumed that the domain
$\Gc_1$ is located on the left of this line. Then one has
\begin{align*}
&\int\left\{f(u)-f(v)-\left[
Q_\lm^u-Q_\lm^v\right]-(u-v)\dot
x_{iii}\right\}dt= \\
&\int\left\{\left(f^\prime(\dots)-\dot x_{iii}\right)\alpha-
\left[Q_\lm^u-Q_\lm^v\right]\right\}dt
\ge\alpha\int\left\{f^\prime(\dots)-\dot x_{iii}\right\}dt\ ,
\end{align*}
because at the level line $u=v+\alpha$ (the domain is located on the left)
$Q_\lm^u=Q(u_x)$, $Q_\lm^v=Q(v_x)$ and $u_x\le
v_x$ but the function $Q$ is monotone.

{\it iv).} Consider the integral with respect to the line $x_{iv}$ of type
d). Without loss of generality again it can be assumed that the domain
$\Gc_1$ is located on the left of this line. Since we have on the line
that the values of $w$ pass from the domain $w>\alpha$ to the domain
$w<\alpha$ as one moves in the positive direction of the axis $x$, we have
the discontinuity with $u^->u^+$ and also $u^-\ge v\ge u^+$. Then
\begin{align*}
&\int\left\{f(u^-)-f(v)-\left[Q_\lm^u-Q_\lm^v
\right]-(u^--v)\dot x_{iv}\right\}dt=\\
&\int\left\{f(u^-)-f(v)+(v-u^-)\dot x_{iv}-
\left[Q_\lm^u-Q_\lm^v\right]\right\}dt
\ge -\int\left[Q^u-Q^v\right]dt
\end{align*}
in consequence of $E$ condition. However at the discontinuity
$u_x\rightarrow -\infty$, therefore $Q_\lm^u\le Q_\lm^v$ whatever is
$Q_\lm^v$.

Thus one has $\sum\limits_k\int\limits_{I_k}(u-v)dx\le
O(\alpha)$. Arguing in a similar way for the domains $w<\alpha$,
ultimately infer that $\int\limits_{\R}|u-v|dx\le O(\alpha)$,
whence as $\alpha\rightarrow 0$ we have $u=v$ a.e. in $\R$ for
every $0<t<T$.

%
%  REMARK 4.1
%
\begin{remark}
The following fact encourages one to obtain that the condition $E$ is
valid for the generalized in the sense of Definition 1 solution (at the
expense of although degenerate but nonzero viscosity). It follows from the
formula (\ref{4}) for the traveling wave solutions that the rise of local
discontinuities with the values $b_{-0}>b_{+0}$
(it is not necessary that $b_{-0}=b_-$ or $b_{+0}=b_+$)
on the left and right of the discontinuity line correspondingly is
possible only providing the following condition is true
$$ f(b)-f(b_{-0})-s(b-b_{-0})\le 0\ \mbox{as}\ b_{+0}\le b\le b_{-0}\ ;\
s=\frac{f(b_{+0})-f(b_{-0})}{b_{+0}-b_{-0}}\ .$$
\end{remark}

{\bf Acknowledgments.} The author thanks A.~S.~Kalashnikov and
A.~Kurga\-nov for useful discussions, and {\it Russian
Foundation for Basic Researches}, grants \No\hspace*{-1.5mm}\No
99-01-00314, 00-01-00387 for financial support.

%
%  REFERENCES
%

%
\vfill\eject
\end{document}